\documentclass[12pt]{amsart}
\usepackage[utf8]{inputenc}
\usepackage{amsmath,amsthm,amssymb,amsfonts,enumerate,color}
\usepackage{amsmath, amsthm, amsfonts, amssymb, color}
\usepackage{mathrsfs}
\usepackage{amsmath, amsthm, amsfonts, amssymb, color}
 \usepackage{mathrsfs}
\usepackage{amsfonts, amsmath}
 \usepackage{amsmath,amstext,amsthm,amssymb,amsxtra}
 \usepackage{txfonts}
 \usepackage[colorlinks, citecolor=blue,pagebackref,hypertexnames=false]{hyperref}
 \allowdisplaybreaks
 \usepackage{pgf,tikz}
\begin{document}

\newcommand{\norm}[1]{\left\Vert#1\right\Vert}
\newcommand{\abs}[1]{\left\vert#1\right\vert}
\newcommand{\set}[1]{\left\{#1\right\}}
\newcommand{\Real}{\mathbb{R}}
\newcommand{\RR}{\mathbb{R}^n}
\newcommand{\supp}{\operatorname{supp}}
\newcommand{\card}{\operatorname{card}}
\renewcommand{\L}{\mathcal{L}}
\renewcommand{\P}{\mathcal{P}}
\newcommand{\T}{\mathcal{T}}
\newcommand{\A}{\mathbb{A}}
\newcommand{\K}{\mathcal{K}}
\renewcommand{\S}{\mathcal{S}}
\newcommand{\blue}[1]{\textcolor{blue}{#1}}
\newcommand{\red}[1]{\textcolor{red}{#1}}
\newcommand{\Id}{\operatorname{I}}

\newtheorem{thm}{Theorem}[section]
\newtheorem{prop}[thm]{Proposition}
\newtheorem{cor}[thm]{Corollary}
\newtheorem{lem}[thm]{Lemma}
\newtheorem{lemma}[thm]{Lemma}
\newtheorem{exams}[thm]{Examples}
\theoremstyle{definition}
\newtheorem{defn}[thm]{Definition}
\newtheorem{rem}[thm]{Remark}

\numberwithin{equation}{section}

\title[Dirichlet problem for Schr\"odinger operators]
{Dirichlet problem for Schr\"odinger operators 
\\ on Heisenberg groups}

 \author[J. Li, Q.Z. Lin, L. Song]{Ji Li, Qingze Lin* and Liang Song}
\thanks{* Corresponding author}

\address{
Ji Li,
Department of Mathematics,
Macquarie University,
Sydney, NSW 2109,
 Australia
}
\email{ji.li@mq.edu.au}

 \address{
   Qingze Lin,
    School of Mathematics,
    Sun Yat-sen University,
    Guangzhou, 510275,
    P.R.~China}
\email{linqz@mail2.sysu.edu.cn}

 \address{
   Liang Song,
    School of Mathematics,
    Sun Yat-sen University,
    Guangzhou, 510275,
    P.R.~China}
\email{songl@mail.sysu.edu.cn}


\begin{abstract}
We  investigate the Dirichlet problem associated to the Schr\"odinger operator $\mathcal L=-\Delta_{\mathbb{H}^n}+V$ on Heisenberg group $\mathbb H^n$:
\begin{align*}
\begin{cases}
\partial_{ss}u(g,s)-\L u(g,s)=0\,,\quad &{\rm in \,\ } \mathbb{H}^n\times\mathbb{R}^+,\\
u(g,0)=f   \,,\quad &{\rm on \,\ } \mathbb{H}^n
\end{cases}
\end{align*}
with $f$ in $L^p(\mathbb{H}^n)$ ($1< p<\infty$) and in $H^1_{\L}(\mathbb{H}^n)$, i.e., the Hardy space associated with $\L$. Here $\Delta_{\mathbb{H}^n}$ is the sub-Laplacian on $\mathbb H^n$ and  the nonnegative potential $V$ belongs to the reverse H\"older class $B_{Q/2}$ with $Q$ the homogeneous dimension of $\mathbb{H}^n$.
The new approach is to
establish a suitable weak maximum principle, which is the key to solve this problem under the condition $V\in B_{Q/2}$. This result is new even back to $\mathbb R^n$ (the condition will become $V\in B_{n/2}$) since the previous known result requires $V\in B_{(n+1)/2}$ which went through a Liouville type theorem.
 \end{abstract}

\maketitle


\section{\bf Introduction and statement of main results}

Let $\Delta=\sum_{i=1}^n\partial_{x_i}^2$ be the  Laplacian on the Euclidean space $\mathbb R^n$. The  harmonic extension of a function is one of the fundamental tools in harmonic analysis and differential equations.  Given a function $f\in L^p(\mathbb R^n)$, a well-known harmonic
extension is the Poisson integral $e^{-t\sqrt{-\Delta}} f=f*p_t$, where $p_t$ is the standard Poisson kernel on $\mathbb R^n$. The following Fatou theorem characterizes the boundary value of the harmonic function $u(x,t)$ on the upper-half space $\mathbb{R}^{n+1}_+$ via the non-tangential maximal function $u^*(x):=\sup\limits_{(y,t): |x-y|<t} |u(y,t)|$ (see for example \cite[Page 119]{St1993}).

\smallskip
\noindent{\bf Theorem A.}{\it\  Suppose $0<p<\infty$.  If $u$ is harmonic in $\mathbb{R}^{n+1}_+$, then $u^*\in L^p(\mathbb R^n)$ if and only if $u$ is the Poisson integral of an $f\in H^p(\mathbb R^n)$, i.e.,
 $$
 u(x,t):= e^{-t\sqrt{-\Delta}} f(x). 
 $$
Moreover, $\|u^*\|_{L^p}\simeq \|f\|_{H^p}$. Here  $H^p (\mathbb R^n)$ is the classical Hardy space.
}

\smallskip

We note that for $1<p<\infty$, $H^p (\mathbb R^n)$ coincides with $L^p (\mathbb R^n)$. For the endpoint $p=\infty$, besides the space $L^\infty (\mathbb R^n)$, the characterization of boundary value in the well-known BMO (bounded mean oscillation) space,
was obtained by Fefferman--Stein \cite{FS72} and Fabes--Johnson--Neri \cite{FJN}.

\smallskip
\noindent{\bf Theorem B.}{\it\     Suppose $u$ is harmonic in $\mathbb{R}^{n+1}_+$. Then $u$ is the Poisson integral of some $f\in {\rm BMO}(\mathbb R^n)$ if and only if $u\in {\rm HMO}$, which means
 $$
 \sup_{x_B\in\mathbb R^n,\ r_B>0} r_B^{-n}\int_0^{r_B}\int_{B(x_B, r_B)}  t|\nabla  u(x,t)|^2 {dx dt }  \leq C.
 $$
 }

The boundary value problems (analogous to Theorems A and B)
have been studied extensively since then. We refer the reader to \cite{FN1,FJN,FN,B13,B15,DYZ,MMMM,JXY,LWZ} and the references therein
 for the further results on this topic.
Recently, Duong--Yan--Zhang \cite{DYZ} extended Theorem B  to the Dirichlet problem for  the Schr\"odinger equation $-\partial_{tt} u+\mathcal{L}u=0$ in $\mathbb{R}_+^{n+1}$ with boundary datum $f\in {\rm BMO}_{\mathcal{L}}$, where $\mathcal{L}=-\Delta+V$, the positive potential $V$ belongs to the reverse H\"older class ${\rm RH}_q$ with $q\geq n$, and ${\rm BMO}_{\mathcal{L}}$ denotes the $\rm BMO$ space associated to $\mathcal{L}$. Later on, the condition $q \geq n$ in \cite{DYZ}  for boundary data $ {\rm BMO}_{\mathcal{L}}$ (i.e., Theorem B)  was improved by Jiang--Li \cite{JL} to $q \geq(n+1) / 2$. The extension of Theorem A was also studied in \cite{DYZ} under the assumption that $q\geq (n+1)/2$ with a modification in the maximal function, i.e., $\|u^*\|_{L^p(\mathbb R^n)}<\infty$ replaced by
$\sup_{t>0}\|u(\cdot,t)\|_{L^p(\mathbb R^n)}<\infty$.

\smallskip

Thus, it is a natural question whether the condition $q\geq (n+1)/2$ in \cite{DYZ}  for boundary data in $L^p(\mathbb R^n), 1<p<\infty,$ or in Hardy space can be improved? Moreover, we also investigate whether such Dirichlet problem in $\mathbb R^{n+1}_+$  holds for sub-Laplacian or more general Schr\"odinger operators in a curved space.

\smallskip

In this paper, we address these questions by studying the Dirichlet problem associated with
the Schr\"odinger  operator
\begin{equation}\label{e1.4}
\L  =-\Delta_{\mathbb{H}^n}+V\,
\end{equation}
on a typical curved space: the Heisenberg group $\mathbb H^n$ (with the homogeneous dimension $Q=2n+2$). We establish an analogue of Theorem A in this setting assuming that the nonnegative  potential $V$ is in the reverse H\"older class with the index $q\geq Q/2$.
This is new even back to $\mathbb R^n$ (which will be $q\geq n/2$), and improves the result in \cite[Prop 2.9]{DYZ} as they required $q\geq (n+1)/2$.

\smallskip

We now state our result in below precisely, and leave the detail notation on the Heisenberg group $\mathbb H^n$ to Section 2.
Recall that a nonnegative locally $L^q\ (1<q<\infty)$ integrable function $V$ on $\mathbb H^n$ is said
to belong to $B_q$ if there exists $C>0$ such that
the reverse H\"older inequality
\begin{eqnarray*}
\left ( \frac{1}{|B|} \int_B V(g)^q \, dg \right )^{\frac{1}{q}}
\leq C \left ( \frac{1}{|B|} \int_B V(g) \, dg \right )
\end{eqnarray*}
holds for every ball $B$ in $\mathbb H^n$.

The Schr\"odinger operators with a potential satisfying the reverse H\"older inequality have been extensively studied. We refer the reader to \cite{Fe,Zh,Shen1994,Shen} for pioneer results, and to \cite{Lu,Li,DZ,LL} for some further developments.


Next, let us recall the definition of $H^p_{\L}(\mathbb{H}^n)$ and $BMO_{\L}(\mathbb{H}^n)$, which were introduced by Lin et al in \cite{LL,LLL}.
For $p\in[1,\infty)$, the Hardy space $H^p_{\L}(\mathbb{H}^n)$ associated to $\L$ is defined by
$$H^p_{\L}(\mathbb{H}^n):=\left\{f\in L^p(\mathbb{H}^n):\, \widetilde{\mathcal{P}}_*(f)\in L^p(\mathbb{H}^n) \right\} \quad {\rm with \ the \ norm} \  \|f\|_{H^p_{\L}(\mathbb{H}^n)}:=\left\|\widetilde{\mathcal{P}}_*(f)\right\|_{L^p(\mathbb{H}^n)}\,. $$
Here $\widetilde{\mathcal{P}}_*(f)$ denotes the nontangential maximal function related to $\sqrt{\L}$, which is defined  by
$$\widetilde{\mathcal{P}}_*(f)(g):=\sup_{d(g,h)<s<\infty}\left|e^{-s\sqrt{\L}}(f)(h)\right|\,,\quad g\in \mathbb{H}^n,$$
where the pseudo-metric $d$ is defined in Section 2.
From the $L^p(\mathbb{H}^n)$  boundedness of the Hardy--Littlewood maximal function ${\mathcal M}$ related to balls in $\mathbb{H}^n$ (see \cite[Chapter~2]{FS}), and the fact  that $\widetilde{\mathcal{P}}_*(f)\leq c{\mathcal M}(f)$, it is not difficult to obtain that $H^p_{\L}(\mathbb{H}^n)=L^p(\mathbb{H}^n)$, for $1<p<\infty$.  It also has been known that the dual space of $H^1_{\L}(\mathbb{H}^n)$ is $BMO_{\L}(\mathbb{H}^n)$\,, which is defined by
$$BMO_{\L}(\mathbb{H}^n):=\left\{f\in L^1_{loc}(\mathbb{H}^n):\, \|f\|_{BMO_{\L}(\mathbb{H}^n)}:=\|f^{\sharp}_V\|_{L^{\infty}(\mathbb{H}^n)}<\infty \right\}\,,$$
where
$$f^{\sharp}_V(g):=\sup_{g\in B:=B(g_B,r_B)}\frac{1}{|B(g_B,r_B)|}\int_{B(g_B,r_B)}\left|f(h)-f_{B,V}\right|dh\,,\quad g\in \mathbb{H}^n\,,$$
and
\begin{equation*}
f_{B,V}:=
\begin{cases}
\ \frac{1}{|B(g_B,r_B)|}\int_{B(g_B,r_B)}f(h)\, dh\,,& {\rm if\ } r_B<\rho(g_B)\,,\\
\ 0\,,&{\rm if\ } r_B\geq\rho(g_B)\,.
\end{cases}
\end{equation*}
Here, $\rho$ is an auxiliary function associated with $V$, defined in Section 2.

For the  theory  of new Hardy and BMO spaces  associated to more general differential operators, we refer to \cite{DZ,DGMTZ,DY1,DY2,AMR,HM,HLMMY} and the references therein.


\smallskip

The nonisotropic Sobolev space $W^{1,2}_{loc}(\mathbb{H}^n\times\mathbb{R}^+)$\,is defined by
$$W^{1,2}_{loc}(\mathbb{H}^n\times\mathbb{R}^+):=\left\{u\in L^{2}_{loc}(\mathbb{H}^n\times\mathbb{R}^+):\, \left|\left(\nabla_{\mathbb{H}^n},\ \partial_s\right) u\right| \in L^{2}_{loc}(\mathbb{H}^n\times\mathbb{R}^+)\right\}\,.$$
Then the equation
$${\mathbb L}u(g,s):=-\partial_s^2u(g,s)+\L u(g,s)=0\,,\quad {\rm in\ } \mathbb{H}^n\times\mathbb{R}^+$$
is understood in the weak sense, that is,
  $u\in {W}^{1, 2}_{{\rm loc}} (\mathbb{H}^n\times\mathbb{R}^+) $ is a weak solution of ${\mathbb L}u =0$   if it satisfies
$$\int_{\mathbb{H}^n\times\mathbb{R}^+} (\nabla_{\mathbb{H}^n},\ \partial_s) u\cdot (\nabla_{\mathbb{H}^n},\ \partial_s)\psi \,dY+
 \int_{\mathbb{H}^n\times\mathbb{R}^+} V u\psi \,dY=0,\ \ \ \ \forall \psi\in C_0^{1}(\mathbb{H}^n\times\mathbb{R}^+)\,.$$
In the sequel, we call such a function $u$ an {\it ${\mathbb L}$-harmonic} function  in the upper-half space $\mathbb{H}^n\times\mathbb{R}^+$.

\smallskip

Let us state our main results, Theorems \ref{the1.5} and \ref{the1.6}, which establishes
an analogue of Theorem A  on Heisenberg group $\mathbb H^n$ for {\it ${\mathbb L}$-harmonic} functions. 

\begin{thm}\label{the1.5}
Suppose $V\in B_{Q/2}$ and $p\in[1,\infty)$\,.  If $f\in H^p_{\L}(\mathbb{H}^n)$\,, then the function $u(g,s):= e^{-s\sqrt{\L}} f(g)$ is the solution of Dirichlet problem
\begin{equation}
\begin{cases}
\ \mathbb{L} u(g,s)=0,&(g,s)\in \mathbb{H}^{n}\times\mathbb{R}^+\,;\\
\ u(g,0)=f(g)\,,&g\in \mathbb H^n;\,\\
\ u^*\in L^p(\mathbb{H}^n),
\end{cases}
\end{equation}
where $u^*(g):=\sup\limits_{d(g,h)<s} |u(h,s)|$, and the boundary data $f$ is taken in the sense of nontangential convergence, i.e.,  $\lim\limits_{\substack{d(g,h)<s\\
s\to 0^+}}u(h,s)=f(g)$ for a.e. $g\in \mathbb{H}^n$.
Moreover, we have $\|u^*\|_{L^p(\mathbb{H}^n)}\leq C \|f\|_{H^p_{\L}(\mathbb{H}^n)}$ and $u(\cdot,s)$ converges to $f$ in $H^p_{\L}(\mathbb{H}^n)$ as $s\to0^+$.
\end{thm}

\begin{thm}\label{the1.6}
Suppose $V\in B_{Q/2}$ and $p\in[1,\infty)$\,. If $u\in C({\mathbb H}^{n}\times\mathbb{R}^+)$ is a weak solution of
${\mathbb L}u =0$ in ${\mathbb H}^{n}\times\mathbb{R}^+$ and $u^*\in L^p(\mathbb{H}^n)$. Then
$u(g,s)= e^{-s\sqrt{\L}} f(g)$ for some $f$ in $H^p_{\L}(\mathbb{H}^n)\,.$
 \end{thm}

\begin{rem}
Comparing to the classical harmonic functions, the ${\mathbb L}$-harmonic functions are no longer smooth in general. And, the  kernel of Poisson semigroup $e^{-t\sqrt{\mathcal{L}}}$  associated to $\L$ will no longer be convolution kernels.   These lead to quite a few  difficulties in the proofs.


It is worth pointing out that Theorems \ref{the1.5} and \ref{the1.6} only require the condition $V\in B_{Q/2}$,
 which is new even in ${\mathbb R}^n$. We refer to \cite[Proposition 2.9]{DYZ}, where $V\in B_{(n+1)/2}$ was required. Let us describe the reason. To show the uniqueness of the solution of ${\mathbb L}$-harmonic equation with Dirichlet boundary value  under  the lack of harmonicity of solutions, Duong et al \cite{DYZ} extended  ${\mathbb L}$-harmonic function from the half-upper space ${\mathbb R}^{n+1}_+$ to the whole space ${\mathbb R}^{n+1}$, then they showed and used a key Liouville type theorem. In the process, the potential $V$ is altered from a function of ${\mathbb R}^{n}$ to that of ${\mathbb R}^{n+1}$. Note that the Liouville type theorem  in \cite{DYZ} required  the condition $V\in B_q$, where $q$ is bigger than a half of dimension. So the condition $V \in B_{(n+1)/2}$ is needed. In \cite{JL}, the authors used the similar argument in the case of  $\mathbb{H}^{n}\times\mathbb{R}^+$. In this paper, we use a new approach to show the uniqueness of the solution of ${\mathbb L}$-harmonic equation with Dirichlet boundary value. We  establish the weak maximum principle on a bounded domain of ${\mathbb H}^{n}\times\mathbb{R}$  for the solution of  $(-\Delta_{\mathbb{H}^n}-\partial_{ss})u=0$   (see Lemma~\ref{WMP}), which, together with a result in \cite{FS}, gives the desired result.
\end{rem}

We also investigate the behavior of the fractional integral operator $\mathcal{I}_\alpha:={\L}^{-\alpha/2}$ with $0<\alpha<Q$\,, which, by functional calculus, has the following expression:
$${\L}^{-\alpha/2}f(g)=\int_0^{\infty} s^{\alpha/2-1}e^{-s\L}f(g) \,ds\,,\quad \forall g\in \mathbb{H}^n\,.$$
For the classical fractional integral operator $I_\alpha:=(-\Delta_{\mathbb{H}^n})^{-\alpha/2}$, it is known (see \cite[Theorem~3.7.2]{SGK}) that $I_\alpha$ is bounded from $L^p(\mathbb{H}^n)$ into $L^q(\mathbb{H}^n)$ for $0<\alpha<Q$ and $1<p<q<\infty$ with $\frac{1}{q}=\frac{1}{p}-\frac{\alpha}{Q}$\,. Therefore, the same conclusion holds for $\mathcal{I}_\alpha$ as well, due to the fact that $0\leq\mathcal{K}_s\leq \mathcal{H}_s$\,, where $\mathcal{K}_s$ is the kernel of the heat semigroup $\big\{ e^{-s\mathcal {L}}:\ s>0
\big\}$ and $\mathcal{H}_s$ is the convolution kernel of
the heat semigroup $\big\{ e^{s\Delta_{\mathbb H^n}}:\ s>0
\big\}$, respectively. Now we consider the endpoint situation: $q=\infty$\,. Just as in the Euclidean space, $I_\alpha$ is not bounded from $L^{Q/\alpha}(\mathbb{H}^n)$ into $L^{\infty}(\mathbb{H}^n)$\,. Instead, the identical proof in \cite[Page~221]{TA} shows that for any $f\in L^{Q/\alpha}(\mathbb{H}^n)$\,, either $I_\alpha f\equiv\infty$ or $I_\alpha f\in {\rm BMO}(\mathbb{H}^n)$. The following theorem will show that for $\mathcal{I}_\alpha$\,, there is only one situation. However, it is not a surprising fact since Dziuba\'{n}ski et al \cite{DGMTZ} has showed the similar result  in ${\mathbb R}^n$.

\begin{thm}\label{fractional}
Suppose $V\in B_{Q/2}.$ If $0<\alpha<Q$\,, then the fractional integral operator $\mathcal{I}_\alpha$ is bounded from $ L^{Q/\alpha}(\mathbb{H}^n)$ into ${\rm BMO}_{\L}(\mathbb{H}^n)$\,.
\end{thm}



\noindent{\bf Concluding remark}: {\it We note that the analogous results in this paper for the more general stratified Lie groups still hold by using the ideas developed here. Moreover, similar results can also be established for the Hardy spaces $H^p_{\L}(\mathbb{H}^n)$ for $0<p<1$, which consist of distributions, by using the ideas in this paper.
To ease the burden of notational complexity, we'll not pursue these parallel results here but leave with the interested readers.}

\medskip

This paper is organized as follows. In Section~\ref{s2}, we recall the basic definitions of Heisenberg group $\mathbb{H}^n$ and several estimates of the heat and Poisson kernels. We show the proof of Theorem~\ref{the1.5} in Section~\ref{s3} and  Theorem~\ref{the1.6} in Section~\ref{s4}. Theorem~\ref{fractional}  will be proven in Section~\ref{s5}.

Throughout the article, the letter ``$C$ " will  denote (possibly different) constants
that are independent of the essential variables. In the context, all functions are assumed to be real-valued, then the same conclusions for complex-valued functions can be deduced by considering their real parts and imaginary parts respectively.

\section{\bf Definitions and Preliminary}\label{s2}

Let's recall the basic definitions of Heisenberg group $\mathbb{H}^n$ (see, for examples, \cite{SGK,St1993,TS}).

The Heisenberg group $\mathbb H^n$ of order $n\in\mathbb{Z}^+$ is a Lie group with underlying manifold $\mathbb R^{2n} \times \mathbb R$ and the group operation $\circ$ is given by
\begin{eqnarray*}
(x, t)\circ(y, s)= \big ( x+y,\, t+s+2 \sum_{j=1}^{n} (x_{n+j} y_{j} -
x_{j} y_{n+j}) \big )\,,
\end{eqnarray*}
where $x,y\in\mathbb R^{2n}$ while $t,s\in \mathbb R$\,.
The canonical basis for the Lie algebra of left-invariant vector fields on
$\mathbb H^n$ is
\begin{eqnarray*}
X_{2n+1}= \frac{\partial}{\partial t},\ \
X_j= \frac{\partial}{\partial x_j}+2x_{n+j}
\frac{\partial}{\partial t},\ \ X_{n+j}=\frac{\partial}{\partial
x_{n+j}}-2x_j \frac{\partial}{\partial t},\ \ j=1, \cdots, n.
\end{eqnarray*}
All non-trivial commutators are $[X_j, X_{n+j}]= -4 X_{2n+1},\, j=1, \cdots, n$.
The sub-Laplacian $\Delta_{\mathbb H^n}$ and the gradient
$\nabla_{\mathbb H^n}$ are defined respectively by
\begin{eqnarray*} \Delta_{\mathbb H^n}= \sum^{2n}_{j=1} X^2_j\qquad
\text{and}\qquad \nabla_{\mathbb H^n}= (X_1, \cdots, X_{2n}).
\end{eqnarray*}
The dilations on $\mathbb H^n$ have the form
\begin{eqnarray*}
\delta_r (x, t)= ( rx, r^2 t), \qquad r>0.
\end{eqnarray*}
The Haar measure on $\mathbb H^n$ coincides with the Lebesgue
measure on $\mathbb R^{2n} \times \mathbb R$. The measure of any
measurable set $E$ is denoted by $|E |$. The homogeneous
norm $\|\cdot\|$ on $\mathbb H^n$ is defined by
\begin{eqnarray*}
\|g\|\, = \big (|x|^4+ 16|t|^2 \big )^\frac{1}{4}, \qquad g= (x,t) \in
\mathbb H^n.
\end{eqnarray*}
This norm satisfies the pseudo-triangle inequality as follows:
\begin{eqnarray}\label{tau}
\|g+h\|\leq \tau (\|g\|+\|h\|),\qquad\forall g,h\in\mathbb H^n,
\end{eqnarray}
where $\tau>1$ is a fixed constant.

Thus, this norm  leads to a
left-invariant distance $$d(g,h) = \|g^{-1}\circ h\|, \qquad\forall g,h\in\mathbb H^n.$$ Then the
ball of radius $r$ centered at $g$ is given by
\begin{eqnarray*}
B(g,r)=\{h \in\mathbb H^n: \; d(g, h )\, <r\}\,,
\end{eqnarray*}
whose volume is given by
\begin{eqnarray*}
\big| B(g,r) \big| =\nu_n r^{Q},\quad \nu_n:= | B(0,1)|\,,
\end{eqnarray*}
where $Q=2n+2$ is the homogeneous dimension of $\mathbb H^n$.

For the given potential $V$\,, the auxiliary function $\rho(g):=\rho(g, V)$ is defined by
\begin{eqnarray*}
\rho(g)= \sup_{r>0}\, \bigg \{ r:\; \frac{1}{r^{Q-2}}\int_{B(g,
r)}V(h)\, dh \leq 1 \bigg \}, \qquad g\in\mathbb H^n\,.
\end{eqnarray*}
The auxiliary function was introduced by Shen in \cite{Shen1994}\,. Basic properties of the auxiliary functions were given
by Shen \cite{Shen1994,Shen} on Euclidean spaces and by Lu \cite{Lu} on
homogeneous spaces. When $V\in B_{Q/2}$, it is known that $0< \rho(g)< \infty$ for any $g \in \mathbb{H}^n$\,.
Moreover, there exists $l_0>0$ such that, for any $g,\,h\in\mathbb H^n$\,,
\begin{eqnarray}\label{auxiliary}
\frac{1}{C} \bigg ( 1+ \frac{d(h,g)}{\rho(g)} \bigg )^{-l_0}
\leq \frac{\rho(h)}{\rho(g)} \leq C \bigg (1+ \frac{d(h,g)}{\rho(g)} \bigg )^{\frac{l_0}{l_0+1}}\,.
\end{eqnarray}
In particular, $\rho(h) \approx \rho(g)$ if $\,d(g,h)< C\,
\rho(g)$\,.

It is well-known that the operator $\L$ generates a ($C_0$) contraction semigroup $\{e^{-t\sqrt{\L}}\}_{t>0}$ on $L^2(\mathbb H^n)$.
Since the potential $V$ is nonnegative, by the Trotter product formula (see \cite[p.~53]{JAG}), the semigroup kernels ${\mathcal K}_s(g,h)$ of the operators $e^{-s\L}$ satisfy
\begin{eqnarray}\label{compareofheat}
0\leq {\mathcal K}_s(g,h)\leq {\mathcal H}_s(h^{-1}\circ g)
\end{eqnarray}
for all $g,h\in \mathbb H^n$ and $s>0$, where $\mathcal{H}_s(g)$ is the convolution kernel of
the heat semigroup $\big\{ e^{s\Delta_{\mathbb H^n}}:\ s>0
\big\}$. The heat kernel $\mathcal{H}_s(g)$ satisfies the estimate
\begin{eqnarray}\label{estimateofheat}
0< \mathcal{H}_s(g) \leq C\, s^{-\frac{Q}{2}} e^{-A_0 s^{-1} \|g\|^2},
\end{eqnarray}
where $A_0$ is a positive constant (see \cite{Jerison} or \cite{LiZhang} for the more precise bounds). An
explicit expression of $\mathcal H_s(g)$ in terms of the Fourier transform
with respect to central variable was given by Hulanicki
\cite{Hulanicki}.

In effect, this estimate for ${\mathcal K}_s(g,h)$ can be improved due to the fact that $V\not\equiv 0$ belongs the reverse H\"older
class $B_{Q/2}$\,. The function $\rho(g)$ arises naturally in this context.

\begin{lem}[\cite{LL,JL}] \label{le2.1} Suppose $V\in B_{Q/2}.$
For every $N>0$ there exist the constants $C_N$ and $c$ such that for $ g,h\in\mathbb H^n, s>0$,
\begin{align*}
&{\rm(1)}\qquad 0\leq {\mathcal K}_s(g,h)\leq C_N s^{-Q/2}e^{-cs^{-1} d(g,h)^2}\left(1+\frac{\sqrt{s}}{\rho(g)}
+\frac{\sqrt{s}}{\rho(h)}\right)^{-N}\,,\\
&{\rm(2)}\qquad \abs{\partial_s{\mathcal K}_s(g,h)}\leq C_Ns^{-\frac{Q+2}{2}}e^{- cs^{-1}d(g,h)^2}
 \left(1+\frac{\sqrt s}{\rho(g)}+\frac{\sqrt s}{\rho(h)}\right)^{-N} .
\end{align*}
\end{lem}

There is a famous Kato--Trotter formula in the perturbation theory for semigroups of operators (see \cite{Kato}):
\begin{align}\label{KTformula}
\mathcal{H}_s(h^{-1}\circ g)- {\mathcal K}_s(g,h)
=  \int_0^s\int_{\mathbb{H}^n} \mathcal{H}_{s-t}(w^{-1}\circ g)V(w){\mathcal K}_{t}(w,h)dwdt.
\end{align}
By formula~\eqref{KTformula}, it follows that
\begin{lem}[see Lemma~9 in \cite{LL}]\label{leHHestimate}
If $V\in B_{Q/2}\,,$ then there exist two positive constants $C,A$ such that
\begin{equation*}
\abs{\mathcal{H}_s(h^{-1}\circ g)- {\mathcal K}_s(g,h)}\leq Cs^{-Q/2}e^{-As^{-1} d(g,h)^2}\cdot\min\left\{\left(\frac{\sqrt{s}}{\rho(g)}\right)^{2-\frac{Q}{q_0}}\,,\,\left(\frac{\sqrt{s}}{\rho(h)}\right)^{2-\frac{Q}{q_0}}\right\},
\end{equation*}
where $g,h\in\mathbb H^n,~s>0$ and $Q/2<q_0<Q\,.$
\end{lem}

The Poisson semigroup associated to $\L$ can be obtained from the heat semigroup through Bochner's
  subordination formula (see \cite{St1970}):
\begin{align}\label{e2.1}
e^{-s\sqrt{\L}}&=\frac{1}{\sqrt{\pi}}\int_0^\infty\frac{e^{-u}}{\sqrt{u}}~e^{-{s^2\over 4u}\L}~du=\frac{s}{2\sqrt\pi}\int_0^\infty
\frac{~e^{-{s^2\over 4t}}}{t^{3/2}}e^{-t\L}~dt.
\end{align}

From \eqref{e2.1},  the Poisson kernels ${\mathcal P}_s(g,h)$, associated to $e^{-s\sqrt{\L}}$,
satisfy the following estimates. For its proof, see \cite[Proposition~3.5]{JL}\,.

\begin{lem}[Proposition~3.5 in \cite{JL}] \label{le2.2} Suppose $V\in B_{Q/2}\,.$

(1)\ For any $N>0$ and $m\in{\{0\}\bigcup\mathbb N}$, there exists a constant $C_{N,\,m}>0$ such that for all $g,h\in\mathbb H^n, s>0$,
\begin{eqnarray}\label{e2.2}
|s^m\partial^m_s{\mathcal P}_s(g,h)|\leq C_{N,\,m} {s \over (s+ d(g,h))^{Q+1}} \left(1+ {s+d(g,h)\over \rho(g)}+{s+d(g,h)\over \rho(h)}\right)^{-N}.
\end{eqnarray}

(2)\ There exists $C>0$ and $q_0\in(Q/2,Q)$ such that
\begin{eqnarray}\label{differenceofPoisson}
|{\mathcal P}_s(g,h)-{\mathcal P}_s(g_0,h)|\leq C \left({ d(g,g_0)\over s+d(g_0,h)}\right)^{2-\frac{Q}{q_0}}{s \over (s+d(g_0,h))^{Q+1}} \,,
\end{eqnarray}
whenever $d(g,g_0)\leq \min\{ d(g_0,h)/4\,, \rho(g_0)\}$ and $s>0$\,.
\end{lem}

\section{\bf Proof of Theorem~\ref{the1.5}}\label{s3}

First, we state the following lemma, which can be found in  \cite[Lemma~8.6]{FS}.
\begin{lem}\label{FSle1}
If $u$ is continuous on $\mathbb H^n\times\mathbb{R}^+$ and $u^*\in L^p(\mathbb H^n)$ for some $p\in [1,\infty)$, then for any $\varepsilon>0$ and $\delta>0$\,, there exist $T>0$ and $R>0$ such that $|u(g,t)|<\varepsilon$ when $t\geq T$ or when $\delta\leq t\leq T$ and $|g|\geq R$\,.
\end{lem}

\medskip

\begin{proof}[Proof of Theorem~\ref{the1.5}]
(1) \ Since $H^p_{\L}(\mathbb H^n)=L^p(\mathbb H^n)$, when $1<p<\infty$ and $H^1_{\L}(\mathbb H^n)\subset L^1(\mathbb H^n)$, we first consider  the case of $f\in L^p(\mathbb H^n)$.

For $f\in L^p(\mathbb H^n)$, $1\leq p<\infty$, we shall show that the Poisson integral $u(g,s):= e^{-s\sqrt{\L}} f(g)$ is a solution of the Dirichlet problem:
\begin{equation*}
\begin{cases}
\ \mathbb{L} u(g,s)=0,&(g,s)\in \mathbb{H}^{n}\times\mathbb{R}^+\,;\\
\ u(g,0)=f(g)\,,&g\in \mathbb H^n\,.
\end{cases}
\end{equation*}

Firstly, we claim  that $e^{-s\sqrt{\L}} f\in W^{1,2}_{loc}(\mathbb{H}^n\times\mathbb{R}^+)$\,. By the estimates \eqref{e2.2} of the Poisson kernel associated to $e^{-s\sqrt{\L}}$ and H\"older's inequality, we know that $e^{-s\sqrt{\L}} f$\,, $\partial_s e^{-s\sqrt{\L}} f$ and $\partial_{ss} e^{-s\sqrt{\L}} f$ are all locally bounded.
For any $\phi\in C^{\infty}_0(\mathbb{H}^n\times\mathbb{R}^+)$\,, we have
\begin{align*}
&\int_{\mathbb{H}^n\times\mathbb{R}^+}|\nabla_{\mathbb{H}^n}e^{-s\sqrt{\L}} f(g)|^2\phi^2 \, dgds\\
&= \int_{\mathbb{H}^n\times\mathbb{R}^+}\left[\left<\nabla_{\mathbb{H}^n}e^{-s\sqrt{\L}} f,\, \nabla_{\mathbb{H}^n}(\phi^2e^{-s\sqrt{\L}}f)\right>-\left<\nabla_{\mathbb{H}^n}e^{-s\sqrt{\L}} f,\, \nabla_{\mathbb{H}^n}(\phi^2)\right> e^{-s\sqrt{\L}} f\right]\,dgds\\
&= \int_{\mathbb{H}^n\times\mathbb{R}^+}\left[\left<\L e^{-s\sqrt{\L}} f,\, \phi^2e^{-s\sqrt{\L}}f\right>-\left<Ve^{-s\sqrt{\L}} f,\, \phi^2e^{-s\sqrt{\L}}f\right>\right]dgds\\
&\hskip 4cm -\frac{1}{2} \int_{\mathbb{H}^n\times\mathbb{R}^+}\left<\nabla_{\mathbb{H}^n}\big((e^{-s\sqrt{\L}} f)^2\big),\, \nabla_{\mathbb{H}^n}(\phi^2)\right>dgds\\
&\leq \int_{\mathbb{H}^n\times\mathbb{R}^+}\left<\partial_{ss} e^{-s\sqrt{\L}} f,\, \phi^2e^{-s\sqrt{\L}}f\right> dg ds+\frac{1}{2} \int_{\mathbb{H}^n\times\mathbb{R}^+}\left(e^{-s\sqrt{\L}} f\right)^2 \left|\Delta_{\mathbb{H}^n}(\phi^2)\right|dgds\\
&<\infty\,.
\end{align*}
Thus, it follows that $e^{-s\sqrt{\L}} f\in W^{1,2}_{loc}(\mathbb{H}^n\times\mathbb{R}^+)$\,.

Secondly, we will prove $u(g,s):= e^{-s\sqrt{\L}} f(g)$ is a weak solution of $\mathbb{L} u(g,s)=0$.
Indeed, for any $\psi\in C^{1}_0(\mathbb{H}^n\times\mathbb{R}^+)$, we have
\begin{align*}
&\int_{\mathbb{H}^n\times\mathbb{R}^+}\left<\nabla_{\mathbb{H}^n}e^{-s\sqrt{\L}} f(g),\, \nabla_{\mathbb{H}^n}\psi(g,s) \right> \, dgds+\int_{\mathbb{H}^n\times\mathbb{R}^+}\left<\partial_s e^{-s\sqrt{\L}} f(g),\, \partial_{s}\psi(g,s) \right> \,dgds\\
&= \int_{\mathbb{H}^n\times\mathbb{R}^+}-\Delta_{\mathbb{H}^n}e^{-s\sqrt{\L}} f(g)\psi(g,s) \, dgds-\int_{\mathbb{H}^n\times\mathbb{R}^+}\partial_{ss} e^{-s\sqrt{\L}} f(g)\psi(g,s) \, dgds\\
&= \int_{\mathbb{H}^n\times\mathbb{R}^+}\left(-\Delta_{\mathbb{H}^n}-\L \right)e^{-s\sqrt{\L}} f(g)\psi(g,s) \, dgds\\
&= \int_{\mathbb{H}^n\times\mathbb{R}^+}-V(g)e^{-s\sqrt{\L}} f(g)\psi(g,s) \, dgds.
\end{align*}

Thirdly, we will prove that for all $f\in L^p(\mathbb H^n)$, $1\leq p<\infty$,
$$\lim\limits_{\substack{d(g,g_0)<s\\
s\to 0^+}}e^{-s\sqrt{\L}} f(g)=f(g_0),  \qquad {\rm for \ \  a.e.} \quad g_0\in \mathbb{H}^n.
$$

To begin with,   we assume that $f\in C^{\infty}_0(\mathbb{H}^n)$. Then for any $g\in \mathbb{H}^n$ satisfying $d(g,g_0)<s<\min\left\{\frac{1}{2\tau}\,, \rho(g_0)\right\}$\, with $\tau$ the constant as in \eqref{tau}, we write
\begin{align*}
e^{-s\sqrt{\L}} f(g)&=e^{-s\sqrt{\L}} (f\chi_{B(g_0,1)^c})(g)+\left(e^{-s\sqrt{\L}}-e^{-s\sqrt{-\Delta_{\mathbb{H}^n}}}\right) (f\chi_{B(g_0,1)})(g)+\left(e^{-s\sqrt{-\Delta_{\mathbb{H}^n}}}\right) (f\chi_{B(g_0,1)})(g)\\
 &=: I_1(g,s)+I_2(g,s)+I_3(g,s)\,.
\end{align*}
For $I_1(g,s)$, by \eqref{e2.2} and H\"older's inequality, we get that
\begin{align*}
\abs{I_1(g,s)}&\leq C\int_{d(g_0,h)\geq1}\frac{s}{d(g,h)^{Q+1}}|f(h)|\, dh \leq C\int_{d(g_0,h)\geq1}\frac{s}{d(g_0,h)^{Q+1}}|f(h)|\, dh\\
&\leq Cs\left(\int_{\mathbb{H}^n}\frac{1}{(1+d(g_0,h))^{(Q+1)p'}} \, dh\right)^{1/p'}\left(\int_{\mathbb{H}^n}|f(h)|^p dh\right)^{1/p} \\
&\leq C_{g_0}s\|f\|_{L^p(\mathbb{H}^n)}\,,
\end{align*}
where $1/p+1/p'=1\,.$ Therefore, $\lim\limits_{\substack{d(g,g_0)<s\\
s\to 0^+}}I_1(g,s)=0\,.$

For $I_2(g,s)$, we note that $d(g,g_0)<s<\rho(g_0)$ implies that $\rho(g)\approx\rho(g_0)$\,.
By Lemma~\ref{leHHestimate} and Bochner's subordination formula \eqref{e2.1}\,,
we see that
\begin{align*}
\abs{I_2(g,s)}&\leq C\frac{s}{2\sqrt\pi}\int_0^\infty
\frac{~e^{-{s^2\over 4t}}}{t^{3/2}} \int_{d(g_0,h)\leq1} \abs{\mathcal{H}_t(h^{-1}\circ g)- {\mathcal K}_t(g,h)}{|f(h)|}~dhdt\\
&\leq C\frac{s}{2\sqrt\pi}\int_0^\infty
\frac{~e^{-{s^2\over 4t}}}{t^{3/2}} \int_{d(g_0,h)\leq1} t^{-Q/2}e^{-A\frac{d(g,h)^2}{t}}\left(\frac{\sqrt{t}}{\rho(g_0)}\right)^{2-Q/q_0}|f(h)|~dhdt\\
&\leq Cs\|f\|_{L^{\infty}(B(g_0,1))}\left(\int_0^{\rho(g_0)^2}+\int_{\rho(g_0)^2}^\infty\right)
\frac{~e^{-{s^2\over 4t}}}{t^{3/2}} \int_{d(g_0,h)\leq1} t^{-Q/2}e^{-A\frac{d(g,h)^2}{t}}\left(\frac{\sqrt{t}}{\rho(g_0)}\right)^{2-Q/q_0}~dhdt\\
&=: J_1+J_2\,.
\end{align*}
For the term $J_1$\,, we take $k=\frac{1}{2}+\varepsilon-\frac{1}{2}\left(2-\frac{Q}{q_0}\right)$ and $\varepsilon\in \left(0,1-\frac{Q}{2q_0}\right)$\,. Then we have
\begin{align*}
J_1&\leq Cs\|f\|_{L^{\infty}(B(g_0,1))}\int_0^{\rho(g_0)^2}
\left(\frac{t}{s^2}\right)^k\frac{1}{t^{3/2}} \left(\frac{\sqrt{t}}{\rho(g_0)}\right)^{2-Q/q_0}dt\\
&\leq C_{g_0}s^{2\left(1-\frac{Q}{2q_0}-\varepsilon\right)}\|f\|_{L^{\infty}(B(g_0,1))}\int_0^{\rho(g_0)^2}
\frac{1}{t^{1-\varepsilon}}dt\\
&\leq C_{g_0}s^{2\left(1-\frac{Q}{2q_0}-\varepsilon\right)}\|f\|_{L^{\infty}(B(g_0,1))}\,.
\end{align*}

For the term $J_2$\,, we have
\begin{align*}
J_2&\leq Cs\|f\|_{L^{\infty}(B(g_0,1))}\int_{\rho(g_0)^2}^\infty
\frac{1}{t^{3/2}} t^{-Q/2}\left(\frac{\sqrt{t}}{\rho(g_0)}\right)^{2-Q/q_0}\int_{d(g_0,h)\leq1} 1~dh\, dt\\
&\leq C_{g_0}s\|f\|_{L^{\infty}(B(g_0,1))}\int_{\rho(g_0)^2}^\infty
\frac{1}{t^{1/2+Q/2+Q/(2q_0)}} dt\\
 &\leq C_{g_0}s\|f\|_{L^{\infty}(B(g_0,1))}.
\end{align*}
Accordingly, we deduce that $\lim\limits_{\substack{d(g,g_0)<s\\
s\to 0^+}}I_2(g,s)=0\,.$

The term $I_3(g,s)$ is a convolution integral. It follows from \cite[Theorem~2.6]{FS} that,  for any $g_0\in \mathbb{H}^n$\,, $\lim\limits_{\substack{d(g,g_0)<s\\
s\to 0^+}}I_3(g,s)=f(g_0)$.
Accordingly, we obtain that  if  $ f\in C^{\infty}_0(\mathbb{H}^n)$, then
\begin{align}\label{Convergence a.e. for C0}
\lim\limits_{\substack{d(g,g_0)<s\\
s\to 0^+}}e^{-s\sqrt{\L}} f(g)=f(g_0), \quad  {\rm for \ any \ g_0\in \mathbb{H}^n.}
\end{align}

In the end,  we extend the above result to all $f\in L^p(\mathbb{H}^n)$\,.  It is easy to see that the nontangential maximal function
$$\widetilde{\mathcal{P}}_*(f)(g_0):=\sup_{d(g,g_0)<s}\left|e^{-s\sqrt{\L}}(f)(g)\right|\,,\quad  g_0\in \mathbb{H}^n\,,$$
is bounded from $L^1(\mathbb{H}^n)$ to $L^{1,\infty}(\mathbb{H}^n)$ and is bounded from $L^p(\mathbb{H}^n)$ to $L^p(\mathbb{H}^n)$ for any $p\in(1,\infty)$, since $\widetilde{\mathcal{P}}_*(f)$ is controlled by Hardy--Littlewood maximal function ${\mathcal M}(f)$.  These facts, together with  \eqref{Convergence a.e. for C0} and a standard  argument, imply that for all $f\in L^p(\mathbb{H}^n)$, it holds that $\lim\limits_{d(g,g_0)<s\to 0^+}e^{-s\sqrt{\L}} f(g)=f(g_0)$ for a.e. $g_0\in \mathbb{H}^n$\,. Since $H^p_{\L}(\mathbb{H}^n)\subset L^p(\mathbb{H}^n)$, $1\leq p<\infty$,
 we have proven the almost everywhere nontangential convergence of $e^{-s\sqrt{\L}}f$ for all $f\in H^p_{\L}(\mathbb{H}^n)$.

\medskip

(2) By the definitions of $H^p_{\L}(\mathbb{H}^n)$ and $u(g,s)$\,, the following inequality naturally holds:
$$\|u^*\|_{L^p(\mathbb{H}^n)}\leq \|f\|_{H^p_{\L}(\mathbb{H}^n)}\,,$$
for $1\leq p<\infty$\,.

Let us turn to consider the convergence in $H^p_{\L}(\mathbb{H}^n)$ norm.
Now we know that $\lim_{s\to 0^+}u(g,s)=f(g)$ for a.e. $g\in \mathbb{H}^n$ and that the nontangential maximal operator $\widetilde{\mathcal{P}}_*$ is bounded from $L^p(\mathbb{H}^n)$ to $L^p(\mathbb{H}^n)$ for any $p\in(1,\infty)$, then by the Lebesgue dominated convergence theorem, it holds that $u(\cdot,s)$ converges to $f$ in $L^p(\mathbb{H}^n)$ as $s\to0^+$ for any $p\in(1,\infty)$\,. Thus, $u(\cdot,s)$ converges to $f$ in $H^p_{\L}(\mathbb{H}^n)$ as $s\to0^+$ for any $p\in(1,\infty)$\,, since $L^p(\mathbb{H}^n)=H^p_{\L}(\mathbb{H}^n)$ when $p\in(1,\infty)$\,.

It remains to show that $f_{s}:=u(\cdot,\, s)= e^{-s\sqrt{\L}} f$ converges to $f$ in $H^1_{\L}(\mathbb H^n)$ as $s\to 0^+$\,. Let $u_s(g,t):=e^{-t\sqrt{\L}} f_s(g)$\,.
It suffices to show that
\begin{align}\label{n.t.m.convergence}
\lim\limits_{s\to 0} \ (u-u_s)^*(g)=0, \ \quad {\rm for \ \ a.e. } \ g\in\mathbb H^n.
\end{align}
In fact, if \eqref{n.t.m.convergence} has been proven, then by Lebesgue's dominated convergence theorem with the domination function $u^*$\,, we deduced that
$$\int_{\mathbb H^n}(u-u_s)^*(g) \,dg\to 0\,,\quad {\rm as\,\,} s\to 0^+\,.$$
That is, $e^{-s\sqrt{\L}} f$ converges to $f$ in $H^1_{\L}(\mathbb H^n)$ as $s\to 0^+$\,.

Now let us prove \eqref{n.t.m.convergence}.
By Lemma~\ref{le2.2}\,, we know that $u$ is continuous in $\mathbb H^n\times\mathbb{R}^+$ whenever $f\in H^1_{\L}(\mathbb H^n)$\,. Then for any $\varepsilon>0$\,, by Lemma~\ref{FSle1}\,, there exists $T>0$ such that when $t\geq T$\,, it holds that $|u(g,t)|<\varepsilon$ for all $g\in\mathbb H^n$\,. Thus,  $|u(g,t)-u_s(g,t)|<2\varepsilon$ for all $g\in\mathbb H^n$ whenever $t\geq T$\,.

Since $\|u^*\|_{L^1(\mathbb H^n)}<\infty$\,, then for almost every $g\in \mathbb H^n$\,, it holds that $u^*(g)<\infty$ and we denote by $\Theta_1$ the set consisting of such points. We denote by $\Theta_2:=\{g\in \mathbb{H}^n : \ \lim\limits_{\substack{d(g,h)<s\\
s\to 0^+}}e^{-s\sqrt{\L}} f(h)=f(g)\}$. Let $\Theta:=\Theta_1\bigcap\Theta_2$. Clearly, $|{\mathbb H^n}\backslash\Theta|=0$.  Now for any fixed $g\in\Theta$\,, we choose $s_1>0$  small enough such that
$$\Bigg|\lim_{\substack{d(g,h)<k\\
k\to0^+}} u(h,k)-u(g,t)\Bigg|<\varepsilon$$
whenever $t\leq 2s_1$\,. Then for any fixed $g\in\Theta$\,, we get that
$$
\sup_{\substack{d(g,h)<t\\
t<s_1}}|u(h,t)-u_{s}(h,t)|\leq 2\varepsilon$$
for any $s<s_1$\,. Now by the uniform continuity of $u$ in the middle bounded closed region $\{(h,t):\, d(g,h)<t\,, s_1\leq t\leq T+s_1 \}$, we can choose $s_2\in (0,s_1)$ small enough such that
$$\sup_{\substack{d(g,h)<t\\
s_1\leq t\leq T}}\left|u(h,t+s)-u(h,t)\right|<\varepsilon$$
whenever $s<s_2$\,. In sum, we have proven for  any $\varepsilon>0$, exists $s_2>0$, such that for all $s<s_2$, there holds
$$\sup_{\substack{d(g,h)<t}}\left|u_s(h,t)-u(h,t)\right|<\varepsilon,
$$
which implies \eqref{n.t.m.convergence}. Thus the proof of Theorem~\ref{the1.5} is complete.
\end{proof}

\smallskip

\section{\bf Proof of Theorem~\ref{the1.6}}\label{s4}

The following  weak maximum principle is key to the proof of Theorem~\ref{the1.6}.
\begin{lem}\label{WMP}
Let $\Omega$ be a bounded domain in $\mathbb{H}^n\times\mathbb{R}$. If $u\in W^{1,2}(\Omega)\bigcap C(\Omega)$ satisfies $$(-\Delta_{\mathbb{H}^n}-\partial_{ss}) u\leq0$$ in $\Omega$ in the weak sense:
\begin{eqnarray*}
\int_{\Omega}(\nabla_{\mathbb{H}^n},\ \partial_s) u(g,s)\cdot (\nabla_{\mathbb{H}^n},\ \partial_s)\varphi(g,s) \, dgds \leq 0\,,\quad \forall {\rm \ nonnegative\ } \varphi\in C^1_0(\Omega),
\end{eqnarray*}
then
$$\sup_{\Omega}u\leq \sup_{\partial\Omega}u\,.$$
\end{lem}

\begin{proof}
For any nonnegative $\varphi\in C^1_0(\Omega)$\,, we have
\begin{eqnarray}\label{eqw1}
\int_{\Omega}(\nabla_{\mathbb{H}^n},\ \partial_s) u(g,s)\cdot (\nabla_{\mathbb{H}^n},\ \partial_s)\varphi(g,s) \, dgds \leq 0\,.
\end{eqnarray}
Suppose, contrary to the assertion, that $\sup_{\Omega}u>\sup_{\partial\Omega}u=u_0\,.$ Then for some constant $C>0$\,, there exists a subdomain $\Omega'\subset\subset\Omega$ in which $v:=u-u_0-C>0$ and $v=0$ on $\partial\Omega'$\,. Then \eqref{eqw1} remains true with $u$ replaced by $v$\,.
Now let $\phi=v$ in $\Omega'$ and $\phi=0$ elsewhere\,. Then $\phi\in W^{1,2}(\Omega)\bigcap C(\Omega)$ and thus can be approximated by functions in $C^1_0(\Omega)$ in the norm of $W^{1,2}(\Omega)$\,. Therefore,
\begin{eqnarray*}
\int_{\Omega'}(\nabla_{\mathbb{H}^n},\ \partial_s) v(g,s)\cdot (\nabla_{\mathbb{H}^n},\ \partial_s) v(g,s) \, dgds \leq 0\,,
\end{eqnarray*}
which implies that $|(\nabla_{\mathbb{H}^n},\ \partial_s)v(g,s)|=0$ for almost every $(g,s)\in\Omega'$\,. By the definition of $\nabla_{\mathbb{H}^n}$\,, it follows that $|(\nabla_{\mathbb{R}^{2n+1}},\ \partial_s)v(g,s)|=0$ for almost every $(g,s)\in\Omega'$\,.
Thus, $v$ is a constant function in $\Omega'$\,. Since $v$ is continuous on $\overline{\Omega'}$ and $v=0$ on $\partial\Omega'$\,, it holds that $v=0$ in $\Omega'$\,, which contradicts the definition of $v$\,. The proof is complete.
\end{proof}

\begin{lem}\label{WMP1}
Let $\Omega$ be a bounded domain in $\mathbb{H}^n\times\mathbb{R}$.
Suppose $ V\in B_{Q/2}\,.$ If $u\in W^{1,2}(\Omega)\bigcap C(\Omega)$ is a weak solution of $\mathbb{L} u=0$ in $\Omega$, then
$$\sup_{\Omega}|u|\leq \sup_{\partial\Omega}|u|.
$$
\end{lem}

\begin{proof}
Since $u\in W^{1,2}(\Omega)\bigcap C(\Omega)$\,, it is clear  that $u^2\in W^{1,2}(\Omega)\bigcap C(\Omega)$\,. For any nonnegative $\varphi\in C^1_0(\Omega)$\,, we have
\begin{align*}
&\int_{\Omega}(\nabla_{\mathbb{H}^n},\ \partial_s) (u^2)\cdot (\nabla_{\mathbb{H}^n},\ \partial_s)\varphi \,dgds \\
&= \int_{\Omega}(\nabla_{\mathbb{H}^n},\ \partial_s) u\cdot (\nabla_{\mathbb{H}^n},\ \partial_s)(2u\varphi) \, dgds -\int_{\Omega}\left[(\nabla_{\mathbb{H}^n},\ \partial_s) u\cdot (\nabla_{\mathbb{H}^n},\ \partial_s) u\right] 2\varphi \, dgds \\
&= -\int_{\Omega}2V u^2\varphi \, dgds- \int_{\Omega}2|(\nabla_{\mathbb{H}^n},\ \partial_s) u|^2 \varphi \, dgds\,\leq\, 0\,.
\end{align*}
Hence, $(-\Delta_{\mathbb{H}^n}-\partial_{ss}) u^2\leq0$. By Lemma \ref{WMP}, we obtain the desired result.
\end{proof}

The following lemma  is a direct consequence of \cite[Lemma~8.7]{FS} and \eqref{e2.2}\,:
\begin{lem}\label{FSle2}
Assume that $f\in L^p(\mathbb H^n)\bigcap C_{00}(\mathbb H^n)$ for some $p\in [1,\infty)$\,, where $C_{00}(\mathbb H^n)$ is the space of continuous functions which vanish at infinity\,. Then for any $\varepsilon>0$\,, there exist $T>0$ and $R>0$ such that $|e^{-t\sqrt{\L}}f(g)|<\varepsilon$ when $t\geq T$ or when $0\leq t\leq T$ and $|g|\geq R$\,.
\end{lem}

We now give the proof of Theorem~\ref{the1.6}.

\begin{proof}[Proof of Theorem~\ref{the1.6}]
Let $\{s_k\}_{k=1}^{\infty}$ be a positive decreasing sequence such that $s_k\to 0$ as $k\to\infty$\,.
Since by assumption $u$ is continuous, then $\lim\limits_{s\rightarrow 0^+} u(h, s+s_k)=u(h, s_k)$ for all $h\in \mathbb{H}^n$.
Since $u^*\in L^p(\mathbb{H}^n)$\,,  by Lemma~\ref{FSle1}\,, we deduce that $u(\cdot,s_k)\in C_{00}(\mathbb{H}^n)$ for all $k\in\mathbb{N}$\,.
Moreover, by the proof of Theorem~\ref{the1.5}, it holds that $\lim\limits_{s\rightarrow 0^+} e^{-s\sqrt{\L}}(u(\cdot, s_k))(h)=u(h, s_k)$ for every $h\in \mathbb{H}^n.$

Define
$$w_k(h,s):=e^{-s\sqrt \L}(u(\cdot, s_k))(h)-u(h,s+s_k).$$
Then $w_k\in C\left(\overline{\mathbb{H}^{n}\times\mathbb{R}^+}\right)$ (due to Lemma~\ref{le2.2} and the assumption that $u$ is continuous and $u^*\in L^p(\mathbb{H}^n)$) and satisfies the following Dirichlet problem:
\begin{equation*}
\begin{cases}
\ \mathbb{L} w_k(g,s)=0,&(g,s)\in \mathbb{H}^{n}\times\mathbb{R}^+\,;\\
\ w_k(g,0)=0\,,&\,g\in \mathbb H^n\,.
\end{cases}
\end{equation*}
Here $\overline{\mathbb{H}^{n}\times\mathbb{R}^+}$ denotes the closure of $\mathbb{H}^{n}\times\mathbb{R}^+$, i.e., $\overline{\mathbb{H}^{n}\times\mathbb{R}^+}=\{(g,t):\ g\in \mathbb{H}^{n},\  t\geq0\}$.

\smallskip
Now we want to show that $w_k$ is identically zero on $\overline{\mathbb{H}^{n}\times\mathbb{R}^+}$\,. Indeed, for any given $\varepsilon>0$\,, since $u(\cdot,s_k)\in C_{00}(\mathbb{H}^n)\bigcap L^p(\mathbb{H}^n)$\,, then by Lemma~\ref{FSle2}\,, there exist $T_1>0$ and $R_1>0$ such that $|e^{-s\sqrt \L}(u(\cdot, s_k))(h)|<\varepsilon$ when $s\geq T_1$ or when $0\leq s\leq T_1$ and $|h|\geq R_1$\,. Applying Lemma~\ref{FSle1} to the function $u$\,, we also get that there exist $T_2>0$ and $R_2>0$ such that $|u(h,s+s_k)|<\varepsilon$ when $s\geq T_2$ or when $0\leq s\leq T_2$ and $|h|\geq R_2$\,.
Take $T:=\max\{T_1,\,T_2\}$ and $R:=\max\{R_1,\,R_2\}$\,. Accordingly, we deduce that for any given $\varepsilon>0$\,, there exist $T>0$ and $R>0$ such that $|w_k(g,s)|<\varepsilon$ when $s\geq T$ or when $0\leq s\leq T$ and $|g|\geq R$\,. In addition, we have known that $w_k(g,0)=0$ for all $g\in \mathbb H^n$\,, therefore, by Lemma~\ref{WMP1}, $|w_k(g,s)|<\varepsilon$ on $\overline{\mathbb{H}^{n}\times\mathbb{R}^+}$\,. Since $\varepsilon$ is arbitrary, it holds that $w_k$ is identically zero on $\overline{\mathbb{H}^{n}\times\mathbb{R}^+}$\,. That is,
\begin{equation}\label{two harmonic function}
u(g,s+s_k)=e^{-s\sqrt \L}(u(\cdot, s_k))(g), \ \ \  (g,s)\in\mathbb{H}^n\times\mathbb{R}^+.
\end{equation}

We first consider the case of $1<p< \infty$\,. Note that as pointed out in  the proof of Theorem~\ref{the1.5}\,, $H^p_{\L}(\mathbb{H}^n)=L^p(\mathbb{H}^n)$ for all $1<p< \infty$\,.
Now since $\|u(\cdot,s_k)\|_{L^p(\mathbb{H}^n)}\leq \|u^*\|_{L^p(\mathbb{H}^n)}\leq C<\infty$
for all $k\in {\mathbb N}$,  there exists a subsequence $\{s_{k_j}\}$ such that $\lim\limits_{j \rightarrow \infty} s_{k_j}=0$, and a function $f\in L^p(\mathbb{H}^n)$
such that $u(\cdot, s_{k_j})$ converges weak* to $f$ as $j\rightarrow \infty.$ That is, for each $\varphi\in L^{p'}(\mathbb{H}^n), 1/p+1/p'=1,$ we have
$$
\lim\limits_{j\rightarrow \infty} \int_{\mathbb{H}^n} u(h, s_{k_j}) \varphi(h)\, dh =\int_{\mathbb{H}^n} f(h) \varphi(h)\, dh.
$$
For any $g\in\mathbb{H}^n, s>0$, we take $\varphi(h)={\mathcal P}_s(g,h)$\,. By Lemma~\ref{le2.2}, $\varphi(\cdot)={\mathcal P}_s(g,\cdot)\in L^{p'}(\mathbb{H}^n)$ for $1\leq p'\leq \infty$\,.
Thus, we have
$$
\lim\limits_{j\rightarrow \infty} \int_{\mathbb{H}^n} {\mathcal P}_s(g,h)u(h, s_{k_j})\, dh =\int_{\mathbb{H}^n} {\mathcal P}_s(g,h) f(h) \,dh.
$$
Therefore, by taking the limits on both sides of \eqref{two harmonic function}, we obtain that for every $g\in\mathbb{H}^n, s>0$
$$
u(g,s)=\lim\limits_{j\rightarrow \infty}
u(g,s+s_{k_j}) =\lim\limits_{j\rightarrow \infty}
e^{-s\sqrt \L}(u(\cdot, s_{k_j}))(g) =e^{-s\sqrt{\L}}f(g).
$$

Now we consider the case of $p=1$\,. Since $\|u(\cdot,s_k)\|_{L^1(\mathbb{H}^n)}\leq \|u^*\|_{L^1(\mathbb{H}^n)}\leq C<\infty$,
for all $k\in {\mathbb N}$,
there exists a finite Borel measure $\mu$ on $\mathbb{H}^n$ which is the weak$\ast$ limit of the sequence $\{ u(\cdot, s_{k_j})\}$ with $s_{k_j}\to0$ as $j\to\infty$. That is,
for each $\varphi$ in $C_{00}(\mathbb{H}^n)$,
$$
\lim\limits_{j\rightarrow \infty} \int_{\mathbb{H}^n} u(h, s_{k_j}) \varphi(h)\, dh =\int_{\mathbb{H}^n} \varphi(h)\, d\mu(h).
$$
By Lemma~\ref{le2.2}, $\varphi(\cdot)={\mathcal P}_s(g,\cdot)$ belongs to $C_{00}(\mathbb{H}^n)$\,.
Thus, we have
$$
\lim\limits_{j\rightarrow \infty} \int_{\mathbb{H}^n} {\mathcal P}_s(g,h) u(h, s_{k_j}) \, dh =\int_{\mathbb{H}^n}  {\mathcal P}_s(g,h) \, d\mu(h)\,.
$$
By taking the limits on both sides of \eqref{two harmonic function}, we have
$$u(g,s)=\int_{\mathbb{H}^n} \mathcal{P}_s(g,h) \, d\mu(h).$$
It remains to prove
$$u(g,s)= e^{-s\sqrt{\L}} f(g)$$ for some $f$ in $H^1_{\L}(\mathbb{H}^n)\,.$
Indeed,
if we can show that $\mu$ is absolutely continuous with respect to Lebesgue measure on $\mathbb H^n$, then $\mu$ coincides with some $f\in L^1(\mathbb H^n)$\,. To see this, we need to prove that for any subset $E\subset\mathbb H^n$\,, the total variation $|\mu|(E)$ of $\mu$ on $E$ is zero whenever $|E|=0$\,.
Since
$$\sup_{s>0}\left|\int_{\mathbb{H}^n} \mathcal{P}_s(\cdot,h) \, d\mu(h)\right|\in L^1(\mathbb{H}^n)\,,$$
for any given $\varepsilon>0$\,, there exists a $\delta>0$ such that for any $F\subset\mathbb H^n$ with $|F|<\delta$\,, it holds that
$$\int_{F}\sup_{s>0}\left|\int_{\mathbb{H}^n} \mathcal{P}_s(g,h) \, d\mu(h)\right| \, dg<\varepsilon\,.$$
Now for any given subset $E\subset\mathbb H^n$ with $|E|=0$\,, there is an open set $\Omega\subset\mathbb H^n$ such that $E\subset\Omega$ and $|\Omega|<\delta$\,. We denote by $C_{00}(\Omega)$ the space of all continuous functions that are supported in $\Omega$ and vanish at infinity. By \cite[(20.49)]{HS}\,, we know that
$$|\mu|(\Omega)=\int_{\Omega}\, d|\mu|=\sup\left\{\left|\int_{\mathbb H^n}\phi \,d\mu\right|:\, \phi\in C_{00}(\Omega)\,, \|\phi\|_{L^{\infty}(\mathbb H^n)}\leq1\right\}\,.$$
Now for every $\phi\in C_{00}(\Omega)$ with $\|\phi\|_{L^{\infty}(\mathbb H^n)} \leq1$\,, by the weak$^*$ convergence in the space of finite Borel measures, we have
\begin{eqnarray*}
\left|\int_{\mathbb H^n}\phi \, d\mu\right| &=&  \lim_{j\to\infty}\left|\int_{\mathbb H^n}\phi(g)\int_{\mathbb{H}^n} \mathcal{P}_{s_{k_j}}(g,h) \, d\mu(h)dg\right|\\
&\leq& \|\phi\|_{L^{\infty}(\mathbb H^n)} \int_{\Omega}\sup_{s>0}\left|\int_{\mathbb{H}^n} \mathcal{P}_s(g,h) \, d\mu(h)\right|dg\\
&<&\varepsilon\,.
\end{eqnarray*}
Thus, $|\mu|(\Omega)<\varepsilon$ and we obtain that the total variation $|\mu|(E)$ of $\mu$ on $E$ is zero whenever $|E|=0$\,, which implies that $\mu$ coincides with some $f\in L^1(\mathbb H^n)$\,, that is, $u(g,s)= e^{-s\sqrt{\L}} f(g)$ for some $f$ in $L^1(\mathbb{H}^n)$\,. Moreover, since
$$
\lim_{j\rightarrow \infty} \int_{\mathbb{H}^n} {\mathcal P}_s(g,h) u(h, s_{k_j}) \, dh =\int_{\mathbb{H}^n}  {\mathcal P}_s(g,h)\, d\mu(h)=\int_{\mathbb{H}^n}  {\mathcal P}_s(g,h)f(h)\, dh,
$$
we obtain that
$$\sup_{\substack{d(g,h)<s\\
s>0}}\left|e^{-s\sqrt{\L}} f(h)\right|\leq \sup_{\substack{d(g,h)<s\\
s>0,\,j>0}}\left|e^{-s\sqrt{\L}} u(h,s_{k_j})\right|= \sup_{\substack{d(g,h)<s\\
s>0,\,j>0}}\left|u(h,s+s_{k_j})\right|\leq u^*(g)\in L^1(\mathbb{H}^n),$$
which shows that $f \in H^1_{\L}(\mathbb{H}^n)$.

Accordingly, we arrive at the conclusion that $u(g,s)= e^{-s\sqrt{\L}} f(g)$ for some $f$ in $H^1_{\L}(\mathbb{H}^n)$\,. The proof is complete.
\end{proof}

\section{\bf Proof of Theorem~\ref{fractional}}\label{s5}

\begin{proof}[Proof of Theorem~\ref{fractional}]
We are to prove that the fractional integral operator $\mathcal{I}_\alpha$ is bounded from $L^{Q/\alpha}(\mathbb{H}^n)$ into ${\rm BMO}_{\L}(\mathbb{H}^n)$\,. The argument is similar to that of \cite[Theorem 7]{DGMTZ} corresponding to the case ${\mathbb R}^n$ with minor modifications.

(1)\ Given any ball $B=B(g_0,r_0)$ with $r_0\geq \rho(g_0)$\,, we need to show that for any $f\in L^{Q/\alpha}(\mathbb{H}^n)$\,, it holds that
$$\frac{1}{|B|}\int_{B}|\mathcal{I}_\alpha f(g)| \, dg\leq C\|f\|_{L^{Q/\alpha}(\mathbb{H}^n)}\,.$$
To this end, we split the function $\mathcal{I}_\alpha f$ as
$$\mathcal{I}_\alpha f(g)=\int_0^{r_0^2} s^{\alpha/2-1}e^{-s\L}f(g)\,ds+\int_{r_0^2}^{\infty} s^{\alpha/2-1}e^{-s\L}f(g)\,ds=:J_1(g)+J_2(g)\,.$$
For the first term $J_1(g)$,  we have
$$|J_1(g)|\leq \int_0^{r_0^2} s^{\alpha/2-1}ds\cdot \sup_{s>0}\left|e^{-s\L}f(g)\right|\leq C_{\alpha}r_0^{\alpha}\cdot {\mathcal M}(f)(g)\,,$$
where ${\mathcal M}$ is the Hardy--Littlewood Maximal operator.
For the second term $J_2(g)$\,, by Lemma~\ref{le2.1} with $N>\alpha$\,, we have
\begin{align*}
|J_2(g)| &\leq C_{\alpha} \int_{r_0^2}^{\infty} s^{\alpha/2-1} \int_{\mathbb{H}^n} s^{-Q/2}e^{-cs^{-1} d(g,h)^2}\left(\frac{\sqrt{s}}{\rho(g)}
\right)^{-N}|f(h)| \, dhds\\
&\leq C_{\alpha}\rho(g)^N {\mathcal M}(f)(g)\int_{r_0^2}^{\infty}s^{\alpha/2-N/2-1}ds\\
&\leq C_{\alpha}r_0^\alpha {\mathcal M}(f)(g),
\end{align*}
where in the last inequality we used the fact that $\rho(g)\leq Cr_0$ if $r_0\geq \rho(g_0)$ and $g\in B(g_0,r_0)$ due to \eqref{auxiliary}.
 Thus, by H\"older's inequality, we deduce that
$$\frac{1}{|B|}\int_{B}|\mathcal{I}_\alpha f(g)| \, dg\leq\frac{C}{|B|}\int_{B}r_0^\alpha {\mathcal M}(f)(g) \,dg \leq\frac{C}{|B|^{1-\alpha/Q}}\int_{B} {\mathcal M}(f)(g) \,dg\leq C\|{\mathcal M}(f)\|_{L^{Q/\alpha}(\mathbb{H}^n)}\,.$$
Then by the boundedness of Hardy--Littlewood Maximal operator ${\mathcal M}$ on $L^{Q/\alpha}(\mathbb{H}^n)$\,, it follows that
$$\frac{1}{|B|}\int_{B}|\mathcal{I}_\alpha f(g)| \, dg\leq C\|f\|_{L^{Q/\alpha}(\mathbb{H}^n)}\,.$$

(2)\ Now we consider the ball $B=B(g_0,r_0)$ with $r_0<\rho(g_0)$\,. Note that in this case, when $g\in B(g_0,r_0)$\,, it holds that $\rho(g)\sim\rho(g_0)$\,. We split the function $\mathcal{I}_\alpha f$ as
$$\mathcal{I}_\alpha f(g)=\int_0^{\rho(g_0)^2} s^{\alpha/2-1}e^{-s\L}f(g)\, ds+\int_{\rho(g_0)^2}^{\infty} s^{\alpha/2-1}e^{-s\L}f(g) \, ds=:J_1(g)+J_2(g)\,.$$
For the second term $J_2(g)$\,, by Lemma~\ref{le2.1} with $N=1$ and H\"older's inequality, we have
\begin{align}\label{J2}
|J_2(g)| &\leq C \int_{\rho(g_0)^2}^{\infty} s^{\alpha/2-1} \int_{\mathbb{H}^n} s^{-Q/2}e^{-cs^{-1}d(g,h)^2}\left(\frac{\sqrt{s}}{\rho(g)}
\right)^{-1}|f(h)| \, dhds\nonumber\\
&= C \rho(g_0)\int_{\rho(g_0)^2}^{\infty} s^{-1/2-1} \int_{\mathbb{H}^n} s^{(\alpha-Q)/2}e^{-cs^{-1}d(g,h)^2}|f(h)| \, dhds\nonumber\\
&\leq C \rho(g_0)\int_{\rho(g_0)^2}^{\infty} s^{-1/2-1} \, ds \cdot \sup_{s>0}\int_{\mathbb{H}^n} s^{(\alpha-Q)/2}e^{-cs^{-1}d(g,h)^2}|f(h)| \, dh\nonumber\\
&\leq C\|f\|_{L^{Q/\alpha}(\mathbb{H}^n)}\,.
\end{align}

It remains to estimate the first term $J_1(g)$\,. To see this, we split $J_{1}(g)$ as
$$J_{1}(g)=\int_0^{r_0^2} s^{\alpha/2-1}e^{-s\L}f(g)\, ds+\int_{r_0^2}^{\rho(g_0)^2} s^{\alpha/2-1}e^{-s\L}f(g) \, ds=:J_{11}(g)+J_{12}(g)\,.$$
For the term $J_{11}(g)$\,, by H\"older's inequality, Minkowski's inequality and \cite[Corollary~2.5]{FS}, we have
\begin{align}\label{J111}
\frac{1}{|B|}\int_{B}|J_{11}(g)| \, dg &\leq \left(\frac{1}{|B|}\int_{B}|J_{11}(g)|^{Q/\alpha} \, dg \right)^{\alpha/Q}\nonumber\\
&\leq \frac{1}{|B|^{\alpha/Q}}\int_0^{r_0^2} s^{\alpha/2-1}ds\cdot\left\|\sup_{s>0}|e^{-s\mathcal{L}}f|\right\|_{L^{Q/\alpha}(\mathbb{H}^n)}\nonumber\\
&\leq C\|f\|_{L^{Q/\alpha}(\mathbb{H}^n)}\,.
\end{align}
For the term $J_{12}(g)$\,, we split the function $f$ into two parts:
$$f=f\chi_{2B}+f\chi_{(2B)^c}=:f_1+f_2\,.$$
Correspondingly, 
$$J_{12}(g)=\int_{r_0^2}^{\rho(g_0)^2} s^{\alpha/2-1}e^{-s\mathcal{L}}f_1(g)\, ds+\int_{r_0^2}^{\rho(g_0)^2} s^{\alpha/2-1}e^{-s\mathcal{L}}f_2(g) \, ds=:J_{121}(g)+J_{122}(g)\,.$$
We now choose $c_B=J_{122}(g_0)$ (which is finite by H\"older's inequality). Then
\begin{align}\label{J112}
\frac{1}{|B|}\int_{B}|J_{12}(g)-c_B| \,dg &\leq \frac{1}{|B|}\int_{B}|J_{121}(g)| \, dg +\frac{1}{|B|}\int_{B}|J_{122}(g)-J_{122}(g_0)| \, dg\,.
\end{align}
For $J_{121}(g)$\,, by Lemma~\ref{le2.1} and H\"older's inequality, we have
\begin{align}\label{J1121}
|J_{121}(g)|&\leq C\int_{r_0^2}^{\rho(g_0)^2} s^{\alpha/2-1}\int_{2B} s^{-\frac{Q}{2}} e^{-A_0 s^{-1} d(g,h)^2}|f(h)|\, dhds\nonumber\\
&\leq C\int_{d(g_0,h)\leq 2r_0}|f(h)| \, dh\int_{r_0^2}^{\infty} s^{\alpha/2-Q/2-1} \, ds\nonumber\\
&\leq C\|f\|_{L^{Q/\alpha}(\mathbb{H}^n)}\,.
\end{align}

It suffices to show that for any $g\in B$,
$$|J_{122}(g)-J_{122}(g_0)|\leq C \|f\|_{L^{Q/\alpha}(\mathbb{H}^n)}\,.$$
In fact, by \cite[Lemma~11]{LL} and H\"older's inequality, there exists a $\delta'\in(0,1)$ such that
\begin{align}\label{J1122}
&|J_{122}(g)-J_{122}(g_0)|\nonumber\\
&\leq C\int_{r_0^2}^{\rho(g_0)^2} s^{\alpha/2-1}\int_{(2B)^c} \left(\frac{d(g,g_0)}{\sqrt{s}}\right)^{\delta'} s^{-\frac{Q}{2}} e^{-A s^{-1} d(g_0,h)^2}|f(h)| \, dhds\nonumber\\
&\leq C r_0^{\delta'}\int_{r_0^2}^{\rho(g_0)^2} s^{-1-\frac{\delta'}{2}}\int_{\mathbb{H}^n} s^{\frac{\alpha-Q}{2}} e^{-A s^{-1} d(g_0,h)^2}|f(h)| \, dhds\nonumber\\
&\leq C r_0^{\delta'}\int_{r_0^2}^{\rho(g_0)^2} s^{-1-\frac{\delta'}{2}}ds\cdot\sup_{s>0}\int_{\mathbb{H}^n} s^{\frac{\alpha-Q}{2}} e^{-A s^{-1} d(g_0,h)^2}|f(h)| \, dh\nonumber\\
&\leq C_N\|f\|_{L^{Q/\alpha}(\mathbb{H}^n)}\,.
\end{align}

Accordingly, by using all estimates from \eqref{J2} to \eqref{J1122}\,, we deduce that for any $B=B(g_0,r_0)$ with $r_0<\rho(g_0)$\,, there is a scalar $c_B$\,, depending on $B$\,, such that
$$\frac{1}{|B|}\int_{B}|\mathcal{I}_\alpha f(g)-c_B|dg\leq C\|f\|_{L^{Q/\alpha}(\mathbb{H}^n)}\,.$$
The proof is complete.
\end{proof}

{\bf Acknowledgements}\ \ The authors are grateful to  Lixin Yan and  Chao Zhang for  helpful suggestions. J. Li is supported by ARC DP 220100285. L. Song and Q.Z. Lin are supported by  NNSF of China (No.~12071490).

 \medskip



\end{document}